\documentclass[11pt]{article}
\usepackage{amsmath,amsthm,amscd,amsfonts,amssymb}
\usepackage{latexsym}
\usepackage[usenames]{color}

\newcommand{\NN}{\mathbb N}

\newcommand{\RR}{\mathbb R}

\newcommand{\ZZ}{\mathbb Z}
\newcommand{\EE}{\mathbb E}

\newcommand{\nn}{\mathcal N}

\newtheorem{theorem}{Theorem}[section]
\newtheorem{lemma}[theorem]{Lemma}
\newtheorem{proposition}[theorem]{Proposition}
\newtheorem{remark}[theorem]{Remark}
\numberwithin{equation}{section}

\newtheorem{conj}[theorem]{Conjecture}
\newtheorem{hyp}[theorem]{Hypothesis}

\def\half{\frac{1}{2}}

\def\pf{{\noindent \bf Proof: }}
\begin{document}
\title{Level Repulsion for a class of decaying random potentials}
\author{Dhriti Ranjan Dolai  and  M Krishna \\
Institute of Mathematical Sciences \\
Taramani, Chennai 600113 \\
dhriti@imsc.res.in, krishna@imsc.res.in \\
}
\maketitle

\begin{center}
{\it Dedicated to Leonid Pastur on his $75^{th}$ birthday }
\end{center}
\newpage
\begin{abstract}
In this paper we consider the Anderson model with decaying randomness and
show that statistics near the band edges in the absolutely continuous 
spectrum in dimensions $d \geq 3$ is independent of the randomness and
agrees with that of the free part.  We also consider the operators
at small coupling and identify the length scales at which the statistics
agrees with the free one in the limit when the coupling constant goes to zero. 
\end{abstract}

\section{Introduction}

In this paper we study the local spectral behaviour in the Anderson
model with decaying randomness in  $d \geq 3$ dimensions.  To do this we 
first consider the non-random part $\Delta$, of the model, restricted to finite
cubes centred at the origin in $\ZZ^d$ and study the eigenvalues of
the resulting matrix.  

The question is motivated by the local statistics of the spectrum
considered for the  Anderson model in the pure point spectral regime.

The studies of statistics of eigenvalues was done in one dimensions by
Molchanov \cite{Mol} and in the Anderson model at high disorder by
Minami \cite{Min} initially.  Both these works formalised the rigorous
procedure for exhibiting Poisson statistics in these random models.  
They show that the eigenvalue statistics near an energy $E$ in the spectrum
follows a Poisson random measure with intensity being $n(E)$ times
the Lebesgue measure, where $n(E)$ is the density of states at $E$. 

Subsequently Poisson statistics was shown for the trees by Aizenman-Warzel 
\cite{AizWar}.  An elegant proof of the Minami estimate needed 
in showing Poisson statistics was obtained by Combes-Germinet-Klein who
also showed Poisson statistics in the continuum models in \cite{ComGerKle}. 
In the paper \cite{GerKlo} Germinet-Klopp give a proof not only of
the Poisson statistics but also showed that the level spacing distribution
follows the exponential law (as one would expect from queueing theory
where the waiting time distribution of a Poissonian queue is exponential).   

In one dimension for a class of decaying random potentials the eigenvalue
statistics was shown to follow the beta-ensemble by Kotani-Nakano \cite{KotNak}. 
Our motivation is to look at the models of decaying random potentials
in $d$ dimension where a sharp mobility edge exists, as shown in 
Kirsch-Krishna-Obermeit \cite{kko} and Jacksic-Last \cite{jl}, and
find out if there is a sharp transition in the local statistics.

It was shown in Dolai \cite{Dhr} that for the models of decaying randomness
in higher dimension which have pure point spectrum, outside $[-2d, 2d]$,
the local statistics is Poisson essentially following the earlier works.
His work combining with this paper shows a transition in the statistics
across the mobility edge.

However in the absolutely continuous spectral regime the statistics is 
different and that is our concern here.  We consider two cases,
one where the random potential is decaying and other where the random potential
at small coupling.   In the former case we identify the rate of decay
and the dimension in which the statistics agrees with that of the 
free operator.  In the latter case we identify the lengths of cubes 
for which the statistics is like the free one.  As far as we know these results 
are new have no comparison in the literature. 

The model we consider is given by 
\begin{equation}\label{themodel}
H^\omega = \Delta + V^\omega, ~ (\Delta u)(n) = \sum_{|m-n|=1} u(m), (V^\omega u)(m) = V^\omega(m) u(m),
\end{equation}
for  $u \in \ell^2(\ZZ^d)$ where $\{V^\omega(n)\}$ is a collection of 
independent real valued random variables. 
We denote the standard basis of $\ell^2(\ZZ^d)$ by
$\{\delta_n, n \in \ZZ^d\}.$  The spectrum $\sigma(\Delta)$ of the operator $\Delta$ is
well known to be purely absolutely continuous and is given by the 
interval $[-2d, 2d]$. We consider 
a cube of side length $2L$ centred at the origin in $\ZZ^d$ 
namely,
$$
\Lambda_L=\{n=(n_1,n_2,\cdots,n_d)\in\mathbb{Z}^d: |n_i|\leq L, i=1,2,\cdots,d \}
$$ 
and take $\delta_{\Lambda_L}$ as the orthogonal projection on to $\ell^2(\Lambda_L)$. 
We define $(2L+1)^d$ dimensional matrices 
$\Delta_L$, $\Delta_{L,E}$ associated with a $E \in (-2d, 2d)$
by
$$
\Delta_L = \chi_{\Lambda_L} \Delta\chi_{\Lambda_L}, ~~ 
\Delta_{L, E} = (L+1)\chi_{\Lambda_L} (\Delta - E) \chi_{\Lambda_L}.
$$
We also consider the matrices
$$
H^\omega_{L, E} = (L+1)\chi_{\Lambda_L} (H^\omega - E) \chi_{\Lambda_L}, ~~ E \in (-2d, 2d).
$$ 
It is known \cite{k1}, \cite{jl}, \cite{kko} that the spectrum of $H^\omega$
is purely absolutely continuous in $(-2d,2d)$ when the variance of $V^\omega(n)$
is finite and the sequence $a_n$ satisfies $a_n \approx |n|^{-2-\epsilon}$ 
as $|n| \rightarrow \infty$.

We then study the measures 
\begin{equation}\label{eqn00}
\mu_{L,E}^0 = \frac{1}{(2L+1)^{d-1}} Tr ( E_{\Delta_{L, E}}() ), ~~~~ 
\mu_{L,E}^\omega = \frac{1}{(2L+1)^{d-1}} Tr ( E_{H^\omega_{L, E}}() )
\end{equation}
where we have notationally denoted the (projection valued ) spectral measure
of a self adjoint operator $A$ by $E_A()$. In the following we also set
$\sigma(A)$ to denote the spectrum of the selfadjoint operator $A$.  

{\bf Acknowledgement:}  We thank Prof S Kotani for his suggestion that we
use martingales and thank Anish Mallik for discussions.

\section{Decaying randomness : Statistics }

In this section we consider perturbations of $-\Delta$ by independent random single site potentials with a either a short range rate of decay at $\infty$
or having a disorder parameter which is small.

\begin{hyp}\label{hyp1}
Let $V^\omega(n) = a_n q^\omega(n)$ with $q^\omega(n)$ independent random variables
distributed according to a probability measure $\nu$ such that 
$\int |x| d\nu(x) < \infty$.  We assume that :
\begin{enumerate}
\item[(i)] the sequence
$a_n$ satisfies $a_n > 0, n \in \ZZ^d$ and $a_n (1+|n|)^{2+\epsilon}$ is bounded.
\item[(ii)] $a_n = \eta, ~ n \in \ZZ^d, ~ \eta >0$.
\end{enumerate}
\end{hyp}

We consider the operators $H^\omega$ as given in the equation (\ref{themodel})
and the measures $\mu_{L,E}^0, \mu_{L, E}^\omega$ given in equation (\ref{eqn00}) associated with the compressions of the operators $\Delta, H^\omega$ to
 $\Lambda_L$.

\begin{theorem}\label{random}
Consider the self adjoint operators $H^\omega$ with $V^\omega$ satisfying 
hypothesis
(\ref{hyp1} (i) ) with the measures
$\mu_{L, E}^\omega$ and $\mu_{L,E}^0$ defined in equation (\ref{eqn00})
associated with $E \in (-2d, 2d)$.  Then for $d \geq 3$  
the sequences of measures $\{\mu_{L,E}^\omega\}$ and $\{\mu_{L,E}^0\}$ have the same 
limit points almost everywhere in the sense of distributions.
\end{theorem}

\pf Note first that positive Radon measure on $\RR$ and positive 
distribution on $C_0^\infty(\RR)$ are the same 
(see Theorem 20.35, \cite{DuiKol}).  
Since $\mu_{L,E}^\omega$, $\mu_{L,E}^0$ are all positive measures,
it is enough to show the convergence in sense of distributions
since the limit points of these then will also be positive distributions
and will be Radon measures. 

For
simplicity we fix $E \in (-2d, 2d)$ and drop the subscript $E$ from
the measures $\mu_{L, E}^\omega, ~ \mu_{L,E}^0$ below.

To this end let $f \in C_0^\infty(\RR)$ and consider
the difference 
$$
\int_{\mathbb{R}} f(x) ~ d\mu^{\omega}_L(x)-\int_{\mathbb{R}} f(x) ~d\mu_L^0(x).
$$
Using the spectral theorem and the definitions of the measures 
$\mu_L^0, \mu_L^\omega$ we can write the above difference as
\begin{equation}\label{eqn10}
\begin{split}
& \int_{\mathbb{R}} f(x) ~ d\mu_L(x)-\int_{\mathbb{R}} f(x) ~d\mu^\omega_L(x)\\
&= \int \widehat{f}(\xi) ~ \frac{1}{(2L+1)^{d-1}}
Tr\left(e^{iH_L^0 \xi} - e^{iH_L^\omega\xi}\right) ~ d\xi \\
\end{split}
\end{equation}
We compute
\begin{equation}\label{eqn11}
\begin{split}
& Tr\left(e^{iH_L^0 \xi} - e^{iH_L^\omega\xi}\right)  \\ 
&= Tr\left(\chi_{\Lambda_L} (e^{iH_L^0 \xi} - e^{iH_L^\omega\xi})\right)  \\ &= 
\int_0^\xi  \sum_{n \in \Lambda_L} \langle \delta_n, e^{iH_L^\omega (\xi - \eta)} i(H_L^\omega - H_L^0) e^{i H_L^0 \eta}\delta_n\rangle ~ d\eta \\
&= \int_0^\xi  \sum_{n, k \in \Lambda_L} \langle \delta_n, e^{iH_L^\omega (\xi - \eta)}\delta_k\rangle  i((L+1) V^\omega(k)\langle \delta_k, e^{i H_L^0 \eta}\delta_n\rangle ~ d\eta \\ 
\end{split}
\end{equation}
Therefore combining the above two equations, we estimate using Cauchy-Schwarz
\begin{equation}\label{eqn12}
\begin{split}
& |\int_{\mathbb{R}} f(x) ~ d\mu_L(x)-\int_{\mathbb{R}} f(x) ~d\mu^\omega_L(x)|\\
& \leq \frac{1}{(2L+1)^{d-1}} \int |(i + \xi) \widehat{f}(\xi)|\\
& ~~\times ~ \frac{1}{|i+\xi|} \int_0^\xi d\eta \sum_{k \in \Lambda_L} (L+1) |V^\omega(k)|  \|e^{i(\xi-\eta) H_L^\omega} \delta_k\| \|e^{i\eta H_L^0 } \delta_k\| \\
& \leq \frac{1}{(2L+1)^{d-2}} \sum_{n \in \Lambda_L} |V^\omega(n)| \int |(i+\xi)\widehat{f}|  ~ d\xi.
\end{split}
\end{equation}

We set
$$
X_L(\omega, f) = \int_{\mathbb{R}} f(x) ~ d\mu^{\omega}_L(x)-\int_{\mathbb{R}} f(x) ~d\mu_L(x).
$$
Then from the above inequality we get the bound
$$
|X_L(\omega, f)|
\leq \|(i+\xi)\widehat{f}\|_1 \frac{1}{(2L+1)^{d-2}}\displaystyle\sum_{n\in\Lambda_L}a_n|q_n(\omega)|.
$$
This estimate and the decay condition on $a_n$ assumed in 
the hypothesis \ref{hyp1} together imply the estimates  
\begin{align}
&  |X_L(\omega,f)| \nonumber
\\ & \leq  \|(i+\xi)\widehat{f}\|_1\frac{1}{(2L+1)^{(d-2)}}\sum_{n\in\Lambda_L}a_n
|q_n(\omega)| \nonumber \\
& \leq C L^{-\frac{\epsilon}{2}} \sum_{n \in \Lambda_L} (1+|n|)^{-d-\frac{\epsilon}{2}} |q^\omega(n)| \nonumber \\
& \leq C L^{-\frac{\epsilon}{2}} \sum_{n \in \Lambda_L} \frac{|q^\omega(n)| - \gamma}{(1+|n|)^{d+\frac{\epsilon}{2}}} + C_k L^{-\frac{\epsilon}{2}} \sum_{n \in \ZZ^d} \frac{\gamma}{(1+|n|)^{d+\frac{\epsilon}{2}}} \label{eqn7} ,
\end{align}
for each fixed $L$ and almost every $\omega$.
We define the random variables 
$$
M_L(\omega) = \sum_{n \in \Lambda_L} (1+|n|)^{-d-\epsilon/2} (|q^\omega(n)| - \gamma), ~~ \mathrm{where} ~~ \gamma = \EE |q^\omega(n)| = \int |x| d\nu(x).
$$
Since $|q^\omega(n)| - \gamma$ are i.i.d  random variables with mean 
zero by hypothesis \ref{hyp1}, we find that the conditional expectation 
of $M_L$ given $M_i, i=1, \dots, L-1$, satisfies
$$
\EE( M_L(\omega) | M_0(\omega), \dots, M_{L-1}(\omega))  =
M_{L-1}(\omega) + \EE(\sum_{|n| = L} (|q^\omega(n)|- \gamma) 
= M_{L-1}(\omega), 
$$
showing that $M_L(\omega)$ is a martingale.  Since
$$
\sup_{L} \EE(M_L(\omega)) < \infty,
$$
the martingale convergence theorem (Theorem 5.7, Varadhan \cite{Var})
shows that $M_L(\omega)$ converges
almost everywhere to a random variable which is finite almost everywhere
which implies that 
$$
L^{-\epsilon/2} M_L(\omega)
$$
converges to zero almost everywhere.  Using this fact in the estimate
(\ref{eqn7}) we find that
$$
|X_L(\omega, f)| 
$$
converges to zero almost everywhere. This estimate is valid for any 
$f \in C_0^\infty(\RR)$, since for functions $f$ in this class
$\|(i+\xi)\widehat{f}\|_1$ is finite.  We define the sequences of 
positive distributions $\Psi_{L, E}^0, \Psi_{L, E}^\omega$ 
$$
\Psi_{L, E}^0(f) = \int f(x) ~ d\mu_{L,E}^0(x), \Psi_{L, E}^\omega(f) = 
\int f(x) ~ d\mu_{L,E}^\omega(x), ~~ f \in C_0^\infty(\RR).
$$
Then from the previous analysis it is clear that $\Psi_{L, E}^0$
and $\Psi_{L, E}^\omega$ have the same limit points almost every $\omega$ as
distributions as desired.  \qed

We now consider the case of weakly coupled random potentials and find the
scales on which the statistics is similar to that of the free part as
the coupling constant goes to zero.  Let $\epsilon(\eta)$ be a function
of $\eta$ such that 
$$
\epsilon(\eta) \rightarrow \infty ~~ \mathrm{if} ~~ \eta \rightarrow 0 ~ \mathrm{and}  ~ 
\lim_{\eta \rightarrow 0} \epsilon(\eta)^2 \eta = 0.
$$

\begin{theorem}\label{aux}
Consider the self adjoint operators $H^\omega$ with $V^\omega$ satisfying
hypothesis
(\ref{hyp1}(ii)), with coupling constant $\eta$.  Consider the measures
$\mu_{L, E}^\omega$ and $\mu_{L,E}^0$ defined in equation (\ref{eqn00})
associated with $E \in (-2d, 2d)$.  Then for $d \geq 1$, 
the sequences of measures $\{\mu_{\epsilon(\eta),E}^\omega\}$ and $\{\mu_{\epsilon(\eta),E}^0\}$ have the same
limit points almost everywhere in the sense of distributions as $\eta \rightarrow 0$.
\end{theorem}
\pf The proof is essentially the same as that of theorem \ref{random}.  In the present case, the first step in the inequality (\ref{eqn7}) becomes, 
\begin{equation}
\begin{split}
|X_{\epsilon(\eta)}(\omega, f)| &\leq \|(1+\xi)\widehat{f}\|_1 \epsilon(\eta)^{-d+2} \sum_{n \Lambda_{\epsilon(\eta)}} \eta |q^\omega(n)| \\
& \leq 
\|(1+\xi)\widehat{f}\|_1 \epsilon(\eta)^2\eta \left( \epsilon(\eta)^{-d} \sum_{n \in \Lambda_{\epsilon(\eta)}} |q^\omega(n)|\right),
\end{split}
\end{equation}
after which the proof is similar to the one given in theorem (\ref{random})
making use of the fact that $\epsilon(\eta)^2\eta \rightarrow 0$ as
$\eta \rightarrow 0$. \qed

\section{\bf Eigenvalues and eigenvectors of $\Delta_{L}$}

In this section we study the eigenvalues of $\Delta_L$ and show that
for energies at the edges of the band $(-2d, 2d)$ there are limit points
for the distributions $\Psi_{L,E}^0$ associated with the measures
$\mu_{L, E}^0$.  

The eigenvalues $\lambda_{j_1, \dots, j_d}^L$ and the (un-normalized) 
eigenfunctions 
$\Psi_{j_1,\dots, j_d, L}$ of $\Delta_L$ are given by
(with the superscript for $\lambda$ denoting an index and not a power )
\begin{equation}\label{evef}
\begin{split}
\lambda_{j_1, \dots, j_d}^L & = 2 \sum_{\ell =1}^d \cos \left(\theta_{j_\ell, L} \right), ~~ \theta_{j, L } = \frac{j \pi}{2(L+1)}, \\ 
\Psi_{j_1, \dots, j_d, L}(n) & = \prod_{\ell=1}^d \phi_{j_\ell, L}(n_\ell), ~ n=(n_1,\dots, n_d) \in \Lambda_L, \\
\phi_{j, L}(m) & = \begin{cases}
\cos\left( \theta_{j, L} m \right), ~~ \mathrm{if} ~~ j ~~\mathrm{is ~ odd}, \\
\sin\left( \theta_{j, L} m \right), ~~ \mathrm{if} ~~ j ~~\mathrm{is ~ even},\end{cases}, ~ m \in \{-L, \dots, L\}, 
\end{split}
\end{equation}
where $j_\ell \in \{1, 2, \dots, 2L+1\}, \ell = 1, \dots, d$.

The eigenvalues of $\Delta_{L, E}$ are correspondingly 
$\{\lambda_{j_1,\dots, j_d}^L - E\}$ for  $E \in [-2d, 2d]$.  

We start with a lemma on the multiplicities of the eigenvalues.
\begin{lemma}\label{lem1}
Let $E_{\Delta_L}$ denote the projection valued measure associated
with $\Delta_L$.  Then for any $\lambda \in \RR$,
$$
Tr(E_{\Delta_L}(\{\lambda\})) \leq d (2L+1)^{d-1}.
$$
\end{lemma}
\pf If $\lambda$ is not an eigenvalue of
$\Delta_L$,  $E_{\Delta_L}(\{\lambda\}) = 0$ and the bound is trivial, 
so we assume without loss of generality that $\lambda \in \sigma(\Delta_L)$. 
The statement in the lemma follows if we show that the
eigenvalues of $\Delta_L$ have multiplicity at most the bound given
in the lemma.  Let  
$$
S = \{ 2 \cos(\frac{k\pi}{2(L+1)}) :  k \in \{1,\dots, 2L+1\}  \}.
$$ 
The points of $S$ are distinct and so $S$ has cardinality $2L+1$ and
the map
$$
f(x_1,\dots, x_d) = x_1 + x_2 + \dots + x_d
$$
from $S^d$ to $[-2d,2d]$ gives precisely all the eigenvalues of $\Delta_L$.  Clearly the
equation $f(x_1, \dots, x_d) = \lambda$ allows the free choice of at most $d-1$
of the variables $x_j$.   If we fix $x_1$ then the number of choices of the remaining
variables is at most $(2L+1)^{d-1}$.  Since we can fix any one of the $d$
variables $x_j$ the bound stated in the lemma follows. \qed 

\begin{remark}\label{rem1}
Since scaling the matrix $\Delta_L$ or adding a constant multiple of the
identity matrix to it does not change the multiplicities of eigenvalues 
the above lemma implies that
$$
Tr( E_{L(\Delta_L - E)}(\{\lambda\}) \leq d (2L+1)^{d-1}.
$$
for any $\lambda \in \RR$.
\end{remark}

\begin{lemma}\label{lem2}
Let $d \geq 1$ and $E \in (-2d, 2d), 2d-2 < |E| < 2d$, then for any 
$f \in C_0^\infty(\RR)$, we have  
$$
\sup_{L \in \NN}  \int f(x)~   d\mu_{L, E}^0(x) < \infty.
$$
\end{lemma}
\pf We give the proof only for the case $2d-2 < E < 2d$, the proof for the
$-2d < E < -2d +2$  is similar. Let $f \in C_0^\infty(\RR)$ have support in $[-K, K]$. 
Let $\Lambda_L^r$ be a cube of side length $L$ in $\ZZ^r$ , take $\Delta_L^0 = 0$
 and set for $r < d$,
$$
(\Delta^r u)(n) = \sum_{|n-i|=1} u(n+i), ~ u \in \ell^2(\ZZ^r), ~~ 
\Delta_L^r = \chi_{\Lambda_{L, r}} \Delta_{r} \chi_{\Lambda_{L,r}}.
$$
Then
\begin{equation}\label{eqn13}
\begin{split}
& \int f(x) d\mu_{L,E}^0(x) \\ &= \frac{1}{(2L+1)^{(d-1)}} \sum_{k=1}^{2L+1} \sum_{\lambda \in \sigma(\Delta_L^{d-1})} f\left((L+1)\left(2\cos(\theta_{k,L}) + \lambda - E\right)\right).
\end{split}
\end{equation} 
The support of $f$ is in $[-K, K]$, so the above sum is only over
$k$ such that $(L+1) \left(2\cos(\theta_{k,L}) + \lambda - E\right)
\in [-K, K]$.  Therefore setting
\begin{equation}\label{eqn13-1}
J_{\lambda, E, L} = [\frac{E-\lambda}{2} - \frac{K}{2(L+1)}, \frac{E -\lambda}{2} + \frac{K}{2(L+1)}], ~ V_{L, r} = (2L +1)^{-r}
\end{equation}
we have
\begin{equation}\label{eqn14}\begin{split}
& |\int f(x) d\mu_{L,E}^0(x)|  \\ & \leq \|f\|_\infty V_{L, d-1} 
\sum_{\lambda \in \sigma(\Delta_L^{d-1})}  \#\left\{k \in \frac{2(L+1)}{\pi} \arccos(J_{\lambda, E, L} \cap [-1, 1])\right\}, 
\end{split}\end{equation}
where 
$$
\arccos(S) = \{ \arccos(x) : x \in S\}.
$$
Letting $|(a, b)| = (b-a)$,noting that the number of integers
in $(a, b)$ is at most $(b-a) + 1$ and using the monotonicity of
$\arccos $ in $[-1, 1]$, the inequality (\ref{eqn14}) becomes
\begin{equation}\label{eqn15}\begin{split}
& |\int f(x) d\mu_{L,E}^0(x)| \\ & \leq  
\|f\|_\infty V_{L,d-1} \sum_{\lambda \in \sigma(\Delta_L^{ d-1}), |\frac{E-\lambda}{2} \pm \frac{K}{2(L+1)} | \leq 1} 1 + \\ & \frac{2(L+1)}{\pi}  \left( \arccos(\frac{E-\lambda}{2} - \frac{K}{2(L+1)}) - \arccos(\frac{E-\lambda}{2} + \frac{K}{2(L+1)}) \right) \\
& \leq \|f\|_\infty +    
\|f\|_\infty V_{L,d-1} \\ & \times \sum_{\lambda \in \sigma(\Delta_L^ {d-1}), |\frac{E-\lambda}{2} \pm \frac{K}{2(L+1)} | \leq 1} \frac{2(L+1)}{\pi}  \left(
\frac{2K}{2(L+1)} \frac{1}{\sqrt{1 - \left(\frac{E-\lambda}{2} + x_L\right)^2}} \right),
\end{split}\end{equation}
where we have used the mean value theorem in the last step for writing
the differences of the arccos terms, which is justified since
$|\frac{E-\lambda}{2}\pm \frac{K}{2(L+1)} | \leq 1$.  By the mean value
theorem it also follows that $|x_L| < \frac{K}{2(L+1)}$ and $ |\frac{E-\lambda}{2} + x_L| < 1$.  
If $d=1$, the proof is over at this stage since for large $L$, the right hand side is  bounded for any $|E| < 1$.  Therefore from now on we assume that $d \geq 2$. Simplifying the above inequality by majorizing it by
the twice the second term, which we can do, otherwise the
proof would be complete,
we get
\begin{equation}\label{eqn16}\begin{split}
& |\int f(x) d\mu_{L,E}^0(x)| \\ & \leq C +   
\|f\|_\infty V_{L,d-1} \sum_{\lambda \in \sigma(\Delta_L^{ d-1}), |\frac{E-\lambda}{2} \pm \frac{K}{2(L+1)} | \leq 1} \frac{2K}{\pi}  \left(
\frac{1}{\sqrt{1 - \left(\frac{E-\lambda}{2} + x_L\right)^2}} \right),
\end{split}\end{equation}
The above term is uniformly bounded in $L$ if 
$(\frac{E-\lambda}{2} + x_L)^2 \leq \half$. So we assume that 
$(\frac{E-\lambda}{2} + x_L)^2 \geq \half$ and in that case 
the sum over $\lambda$ splits into two parts, according as
$\pm(\frac{E-\lambda}{2} + x_L) > \half$.  
Therefore we set
$$
I_{\pm} = 
\|f\|_\infty V_{L,d-1} \sum_{\lambda \in \sigma(\Delta_L^{ d-1}), \pm\frac{E-\lambda}{2} \pm \frac{K}{2(L+1)} \geq \half} \frac{2K}{\pi}  \left(
\frac{1}{\sqrt{1 - \left(\frac{E-\lambda}{2} + x_L\right)}} \right).
$$
We continue with the proof for $I_+$  the proof of the other case
is similar. We have 
$$
\frac{1}{\sqrt{1 - \left(\frac{E-\lambda}{2} + x_L\right)^2}} \leq 
\frac{1}{\sqrt{ 1 - \frac{E-\lambda}{2} - x_L}\sqrt{ 1 + \frac{E-\lambda}{2} + x_L}}
\leq \frac{1}{\sqrt{ 1 - \frac{E-\lambda}{2} - x_L}}.
$$
Using this bound we find 
\begin{equation}\label{eqn17}\begin{split}
& I_+ \leq
\|f\|_\infty V_{L,d-1} \sum_{\lambda \in \sigma(\Delta_L^{ d-1})} \frac{2K}{\pi}  \left(
\frac{1}{\sqrt{1 - \frac{E-\lambda}{2} - x_L}} \right),
\end{split}\end{equation}
We now use the fact that $\lambda \in \sigma(\Delta_L^{ d-1})$ can be split
into $\lambda = \lambda_1 + \lambda_2 $, where
$\lambda_2 \in \sigma(\Delta_L^1)$ and
$\lambda_1 \in \sigma(\Delta_L^{d-2})$. Then the above inequality becomes
\begin{equation}\label{eqn18}\begin{split}
& I_+ \leq
\frac{2K\|f\|_\infty}{\pi} V_{L,d-2} \sum_{\lambda_1 \in \sigma(\Delta_L^{d-2})}  \\ & \times \frac{1}{(2L+1)} \sum_{\lambda_2 \in \sigma(\Delta_L^1), E-\lambda_1 -\lambda_2 -2 x_L > 0}  \left(
\frac{1}{\sqrt{1 - \frac{E-\lambda_1 - \lambda_2 }{2} - x_L}} \right),
\end{split}\end{equation} 
We claim that the
sum 
$$
I(\gamma) =  \frac{1}{(2L+1)} \sum_{\lambda_2 \in \sigma(\Delta_L^1), \lambda_2 < 2\gamma}  \left(
\frac{1}{\sqrt{\gamma + \frac{ \lambda_2 }{2}}} \right)
$$
where $\gamma = 1 - \frac{E -\lambda_1}{2}$, is uniformly bounded in $\gamma$ and $L$. If the claim is true
then we get the bound
\begin{equation}\label{eqn20}\begin{split}
& I_+ \leq
\frac{2K\|f\|_\infty}{\pi} V_{L,d-2} \sum_{\lambda_1\in \sigma(\Delta_L^{d-2})} C <  \frac{2K \|f\|_\infty}{\pi} C, 
\end{split}\end{equation}
giving the lemma.  We therefore prove the claim.  
Using the explicit expressions for the points in $\sigma(\Delta_L^1)$ we 
computed earlier in equation (\ref{evef}), we get
$$
I(\gamma) = \frac{1}{(2L+1)} \sum_{k= 1, \gamma > \cos(\frac{k\pi}{2(L+1)})}^{2L+1}  \left(
\frac{1}{\sqrt{\gamma - \cos\frac{k \pi}{2(L+1)}  }} \right)
$$
Since the function 
$$
g(x) = \frac{1}{\sqrt{\gamma - \cos(x \pi) }}
$$
is monotonically decreasing in  $0 \leq x \leq 1$  we bound the sum above 
by the integral
$$
I(\gamma) \leq \delta_L + \left(\frac{2(L+1)}{2L+1}\right) \int_0^1 
g(x)\chi_{(([-1,\gamma))}(\cos(x\pi)) ~ dx.  
$$
where $\delta_L$ is a small error that is uniformly bounded in $L$.
Changing variables $y= \cos(x\pi)$ gives the bound
\begin{equation}\label{eqn19}\begin{split}
I(\gamma) &\leq \delta_L + \frac{2}{\pi} \int_{-1}^\gamma  g(y)  \frac{1}{\sqrt{1 - y^2}} dy \\
& \leq \delta_L + \frac{2}{\pi} \int_{-1}^\gamma \frac{1}{\sqrt{\gamma - y}} \frac{1}{\sqrt{1-y^2}} dy \\
& \leq \delta_L + \frac{1}{\sqrt{1 - \gamma}} \int_{-1}^\gamma \frac{1}{\sqrt{\gamma - y}}\frac{1}{\sqrt{1+y}}.
\end{split}\end{equation}
The condition on $E$ assures us that $\gamma <0$, therefore the factor
$\frac{1}{\sqrt{1-\gamma}}$ is bounded by $1$, on the other hand a bound by
splitting the integral into two pieces up to and from the midpoint $(1 +|\gamma|)/2$ yields 
$$
 \int_{|\gamma|}^1 \frac{1}{\sqrt{ (y - |\gamma|)(1 - y)} }dy \leq 2. 
$$
Note: In case $|\gamma | =1$ we define $I(\gamma)$ to be
$$
I(\gamma) =  \lim_{\epsilon \downarrow } I(\gamma, \epsilon) 
$$
where
$$
I(\gamma, \epsilon) =  \frac{1}{(2L+1)} \sum_{\lambda_2 \in \sigma(\Delta_L^1), |\lambda_2 | < 2 - \epsilon}  \left(
\frac{1}{\sqrt{\gamma + \frac{ \lambda_2 }{2}}} \right)
$$ 
and bound  $I(\gamma, \epsilon)$ for each $\epsilon >0$, which we can do since all the terms are finite for each $\epsilon >0$.  This avoids the logarithmic singularity in the integral when we
replace the sum defining $I(\gamma)$ by an integral. 
 \qed

\begin{proposition}\label{prop1}
The measures $\mu_{L,E}^0$ have limit points in the vague sense
when $2d-2 < |E| < 2d$. 
\end{proposition}
\pf By lemma above the measures $\mu_{L, E}^0$ are uniformly bounded on the
space of continuous functions of compact support, hence by Helly's selection
theorem they have limit points in the vague sense (by a diagonal argument
if necessary).   To show that there is at least one non-trivial limit
point we show that for some positive function of compact support,
$$
\liminf_{L \in \NN}  f(x) d\mu_{L,E}^0(x) > 0.
$$
To this end consider a $K > 1$ fixed and let $0 \leq f\leq 1$ be 
a continuous function with $f(x) = 1, -K \leq x \leq K$. 
Then we see from equations (\ref{eqn13}), (\ref{eqn13-1}) that 
$$
\int f(x) d\mu_{L,E}^0(x) \geq 
\sum_{\lambda \in \sigma(\Delta_L^{d-1})}  \#\left\{k \in \frac{2(L+1)}{\pi} \arccos(J_{\lambda, E, L} \cap [-1, 1])\right\}. 
$$
As estimated in equation (\ref{eqn14}) we estimate the number of integers $k$
by the Lebesgue measure of the interval, but now taking a smaller interval
$[\frac{E-\lambda}{2} - \frac{K}{4(L+1)}, \frac{E-\lambda}{2} + \frac{K}{2(L+1)}]$ to get as in equation (\ref{eqn15}) (now for lower bound)
\begin{equation}\begin{split}\label{eqn21}
& \int f(x) d\mu_{L,E}^0(x) \\ & \geq 
\sum_{\lambda \in \sigma(\Delta_L^{d-1})}   
\frac{2(L+1)}{\pi}  \left( \arccos(\frac{E-\lambda}{2} - \frac{K}{4(L+1)}) - \arccos(\frac{E-\lambda}{2} + \frac{K}{4(L+1)}) \right) 
\end{split}\end{equation}
using the monotonicity of $\arccos$ in $[-1,1]$. For some $0 < \delta < 1$,  
we take $L$ large so that $\frac{K}{4(L+1)} < \delta/4$, hence using the mean value theorem we get the lower bound 
\begin{equation}\begin{split}
& \frac{2(L+1)}{\pi}  \left( \arccos(\frac{E-\lambda}{2} - \frac{K}{4(L+1)}) - \arccos(\frac{E-\lambda}{2} + \frac{K}{4(L+1)}) \right) 
\\ & = \frac{2(L+1)}{\pi} \frac{K}{2(L+1)} \frac{1}{\sqrt{1 - \left(\frac{E-\lambda}{2} + x_L\right)^2}} \geq \frac{K}{\pi} 
\frac{1}{\sqrt{1 - \left(\frac{E-\lambda}{2} + x_L\right)^2}}. 
\end{split}\end{equation}
Therefore from equation (\ref{eqn21}) and the above we get since
$|x_L| < \delta/4$ for large enough $L$,
$$
 \int f(x) d\mu_{L,E}^0(x)  \geq 
\frac{K}{\pi (2L+1)^{d-1}}\sum_{\lambda \in \sigma(\Delta_L^{d-1}) ~~    
\frac{|E - \lambda|}{2} \leq 1 -  \delta/2 } \sqrt{\frac{1}{2}}.
$$
The right hand side clearly has a limit in terms of the density of states
of $\Delta^{d-1}$ namely
$$
\frac{K}{\pi \sqrt{2}} \nn_{d-1}((E - 2 + \delta, E + 2 - \delta))
$$
where $\nn_{r}$ is the density of states of $\Delta^r$.  For 
$|E| \in (2d-2, 2d)$, 
$(E -2 +\delta, E+2 - \delta) \cap (-2d+2, 2d-2) \neq \emptyset$ for 
small enough $\delta$ showing the positivity of the right hand side.
 \qed

We end this paper with a 
\begin{conj}
If $d \geq 4$, the limit points of $\mu_{L,E}^0$ are independent of 
$E \in (-2d, 2d)$ and they are given by
$$
\sum_{k \in \ZZ} \int \sin(\theta) n_{d-1}(E -2\cos(\theta)) \delta_{\pi k \sin(\theta)}\left( \right) d\theta, 
$$
where $n_{d}$ is density of states of $\Delta$ in $d$ dimensions.  
\end{conj}

We note that the density of states of $\Delta$ is the density of the 
measure $\langle \delta_0, E_{\Delta}() \delta_0\rangle$ and in 
$d\geq 4$ this density function is continuously differentiable, since
its Fourier transform $\widehat{n_{d}}(t)$ is bounded and decays like 
$|t|^{-d/2}$ as can be seen by putting together Lemma 4.1.8, 4.1.9 \cite{DemKri}
via the spectral theorem.

\end{document}